\documentclass[10pt,reqno]{smfart}
\usepackage{amsmath,amssymb,smfthm,stmaryrd,epic}
\usepackage[frenchb]{babel}
\usepackage[latin1]{inputenc}
\usepackage[all]{xy}
\usepackage{a4wide}

\usepackage{xypic}
\usepackage[active]{srcltx}
\usepackage{rotating}

\NumberTheoremsIn{subsection}\SwapTheoremNumbers

\begin{document}

\newtheorem{prop-defi}[smfthm]{Proposition-Définition}
\newtheorem{theo-defi}[smfthm]{Théorème-définition}
\newtheorem{lem-defi}[smfthm]{Lemme-définition}
\newtheorem{notas}[smfthm]{Notations}
\newtheorem{nota}[smfthm]{Notation}
\newtheorem{defis}[smfthm]{Définitions}
\newtheorem{remas}[smfthm]{Remarques}

\newtheorem{theob}{Théorème}[section]
\def\thetheob{\arabic{section}.\arabic{theob}}
\newtheorem{propb}[theob]{Proposition}
\newtheorem{lemb}[theob]{Lemme}
\newtheorem{corob}[theob]{Corollaire}
\newtheorem{defib}[theob]{Définition}
\newtheorem{defisb}[theob]{Définitions}
\newtheorem{remab}[theob]{Remarque}

\renewcommand{\theequation}{\Roman{part}.\arabic{section}.\arabic{subsection}.\arabic{equation}}

\def\Am{{\mathbb A}}
\def\Fm{{\mathbb F}}
\def\Mm{{\mathbb M}}
\def\Nm{{\mathbb N}}
\def\Pm{{\mathbb P}}
\def\Qm{{\mathbb Q}}
\def\Zm{{\mathbb Z}}
\def\Dm{{\mathbb D}}
\def\Cm{{\mathbb C}}
\def\Rm{{\mathbb R}}
\def\Gm{{\mathbb G}}

\def\AC{{\mathcal A}}
\def\CC{{\mathcal C}}
\def\DC{{\mathcal D}}
\def\EC{{\mathcal E}}
\def\FC{{\mathcal F}}
\def\GC{{\mathcal G}}
\def\HC{{\mathcal H}}
\def\IC{{\mathcal I}}
\def\JC{{\mathcal J}}
\def\KC{{\mathcal K}}
\def\LC{{\mathcal L}}
\def\MC{{\mathcal M}}
\def\NC{{\mathcal N}}
\def\OC{{\mathcal O}}
\def\PC{{\mathcal P}}
\def\UC{{\mathcal U}}
\def\VC{{\mathcal V}}
\def\XC{{\mathcal X}}
\def\YC{{\mathcal Y}}

\def\AF{{\mathfrak A}}
\def\GF{{\mathfrak G}}
\def\EF{{\mathfrak E}}
\def\CF{{\mathfrak C}}
\def\DF{{\mathfrak D}}
\def\JF{{\mathfrak J}}
\def\LF{{\mathfrak L}}
\def\MF{{\mathfrak M}}
\def\NF{{\mathfrak N}}
\def\XF{{\mathfrak X}}
\def\UF{{\mathfrak U}}
\def\KF{{\mathfrak K}}

\def \longmapright#1{\smash{\mathop{\longrightarrow}\limits^{#1}}}
\def \mapright#1{\smash{\mathop{\rightarrow}\limits^{#1}}}
\def \lexp#1#2{\kern \scriptspace \vphantom{#2}^{#1}\kern-\scriptspace#2}
\def \linf#1#2{\kern \scriptspace \vphantom{#2}_{#1}\kern-\scriptspace#2}
\def \linexp#1#2#3 {\kern \scriptspace{#3}_{#1}^{#2} \kern-\scriptspace #3}

\def \a{\alpha}
\def \b{\beta}
\def \d{\delta}
\def \e{\epsilon}
\def \g{\gamma}
\def \k{\kappa}
\def \l{\lambda}
\def \m{\mu}
\def \n{\nu}
\def \o{\omega}
\def \r{\rho}
\def \s{\sigma}
\def \t{\tau}
\def \th{\theta}
\def \u {\upsilon}
\def \x{\chi}
\def \vphi {\varphi}

\let \leq=\leqslant
\let \geq=\geqslant
\def \lefto{\longleftarrow}
\def \fin{\hfill $\square$}
\let \DS=\displaystyle
\let \SS=\scriptstyle
\let \longto=\longrightarrow
\let \oo=\infty

\def \FH{\mathop{\mathrm{FH}}\nolimits}
\def \FPH{\mathop{\mathrm{FPH}}\nolimits}
\def \coh{\mathop{\mathrm{Coh}}\nolimits}
\def \res{\mathop{\mathrm{res}}\nolimits}
\def \op{\mathop{\mathrm{op}}\nolimits}
\def \rec {\mathop{\mathrm{rec}}\nolimits}
\def \art{\mathop{\mathrm{Art}}\nolimits}
\def \hyp {\mathop{\mathrm{Hyp}}\nolimits}
\def \cusp {\mathop{\mathrm{Cusp}}\nolimits}
\def \Iw {\mathop{\mathrm{Iw}}\nolimits}
\def \JL {\mathop{\mathrm{JL}}\nolimits}
\def \speh {\mathop{\mathrm{Speh}}\nolimits}
\def \isom {\mathop{\mathrm{Isom}}\nolimits}
\def \Vect {\mathop{\mathrm{Vect}}\nolimits}
\def \groth {\mathop{\mathrm{Groth}}\nolimits}
\def \lef {\mathop{\mathrm{Lef}}\nolimits}
\def \fix {\mathop{\mathrm{Fix}}\nolimits}
\def \hom {\mathop{\mathrm{Hom}}\nolimits}
\def \deg {\mathop{\mathrm{deg}}\nolimits}
\def \val {\mathop{\mathrm{val}}\nolimits}
\def \det {\mathop{\mathrm{det}}\nolimits}
\def \rep {\mathop{\mathrm{Rep}}\nolimits}
\def \spec {\mathop{\mathrm{Spec}}\nolimits}
\def \fr {\mathop{\mathrm{Fr}}\nolimits}
\def \frob {\mathop{\mathrm{Frob}}\nolimits}
\def \ker {\mathop{\mathrm{Ker}}\nolimits}
\def \im {\mathop{\mathrm{Im}}\nolimits}
\def \Red {\mathop{\mathrm{Red}}\nolimits}
\def \red {\mathop{\mathrm{red}}\nolimits}
\def \aut {\mathop{\mathrm{Aut}}\nolimits}
\def \diag {\mathop{\mathrm{diag}}\nolimits}
\def \spf {\mathop{\mathrm{Spf}}\nolimits}
\def \Def {\mathop{\mathrm{Def}}\nolimits}
\def \twist {\mathop{\mathrm{Twist}}\nolimits}
\def \supp {\mathop{\mathrm{Supp}}\nolimits}
\def \Id {{\mathop{\mathrm{Id}}\nolimits}}
\def \bar {\overline}
\def \ind {\mathop{\mathrm{Ind}}\nolimits}
\def \mod {\mathop{\mathrm{mod}}\nolimits}
\def \ker {\mathop{\mathrm{Ker}}\nolimits}
\def \coker {\mathop{\mathrm{Coker}}\nolimits}
\def \mult {\mathop{\mathrm{mult}}\nolimits}
\def \vide{\emptyset}
\def \bad {{\mathop{\mathrm{Bad}}\nolimits}}
\def \gal {{\mathop{\mathrm{Gal}}\nolimits}}
\def \Nr {{\mathop{\mathrm{Nr}}\nolimits}}
\def \rn {{\mathop{\mathrm{rn}}\nolimits}}
\def \vol {{\mathop{\mathrm{vol}}\nolimits}}
\def \ad {{\mathop{\mathrm{ad}}\nolimits}}
\def \tr {{\mathop{\mathrm{Tr~}}\nolimits}}
\def \Sp {{\mathop{\mathrm{Sp}}\nolimits}}
\def \lie {{\mathop{\mathrm{Lie}}\nolimits}}
\def \st {{\mathop{\mathrm{Sp}}\nolimits}}
\def \sp{{\mathop{\mathrm{Sp}}\nolimits}}
\def \card{{\mathop{\mathrm{card}}\nolimits}}
\def \sym{{\mathop{\mathrm{Sym}}\nolimits}}
\def \perv{\mathop{\mathrm{Perv}}\nolimits}
\def \sh {{\mathop{\mathrm{Sh}}\nolimits}}
\def \const {{\mathop{\mathrm{Const}}\nolimits}}
\def \Ind {{\mathop{\mathrm{Ind}}\nolimits}}

\def \ele{élément }
\def \eles{éléments }
\def \cad{c'est à dire }
\def \rem{{\noindent\textit{Remarque:~}}}
\def \exem{{\noindent \textit{Exemple:~}}}
\def \ssi{~si et seulement si~}
\def \cl {{\mathop{\mathrm{cl}}\nolimits}}
\def \Tw {{\mathop{\mathrm{Tw}}\nolimits}}
\def \ob {{\mathop{\mathrm{Ob}}\nolimits}}
\def \ext {{\mathop{\mathrm{Ext}}\nolimits}}
\def \End {{\mathop{\mathrm{End}}\nolimits}}
\def \inv {{\mathop{\mathrm{inv}}\nolimits}}
\def \fix {{\mathop{\mathrm{Fix}}\nolimits}}

\def\semi{\mathrel{>\!\!\!\triangleleft}}


\setcounter{secnumdepth}{3} \setcounter{tocdepth}{3}

\newcommand{\marque}{\addtocounter{smfthm}{1}
{\smallskip \noindent \textit{\thesmfthm}~---~}}

\renewcommand\atop[2]{\ensuremath{\genfrac..{0pt}{1}{#1}{#2}}}

\title{Sur l'irréductibilité de quelques variétés d'Igusa}

\alttitle{On irreductibility of somme Igusa varieties}

\author{Boyer Pascal}
\email{boyer@math.jussieu.fr}
\address{Institut de mathématiques de Jussieu \\ UMR 7586, université Paris 6 \\
175 rue du Chevaleret Paris 13}

\frontmatter

\begin{abstract} On prouve que les variétés d'Igusa de seconde espèce
définies par Harris et Taylor dans leur livre \cite{h-t}, sont aussi connexes que possible, i.e.
que les composantes connexes sont naturellement en bijection avec les caractères du groupe des
inversibles de l'ordre maximal de l'algèbre à division centrale qui lui est attachée.
\end{abstract}

\begin{altabstract}
We describe the connected components of Igusa's varieties of second species defined par Harris and
Taylor in their book \cite{h-t}. There are in bijection with the characters of the inversible
group of the maximal order of the division algebra attached to them.
\end{altabstract}



\altkeywords{Shimura varieties, formal modules, Langlands correspondences, Jacquet-Langlands correspondences,
monodromy filtration, weight-monodromy conjecture, perverse sheaves, vanishing cycles}

\maketitle

\tableofcontents

\pagestyle{headings} \pagenumbering{arabic}

\section*{Introduction}
\renewcommand{\theequation}{\arabic{equation}}
\backmatter

A l'origine de cette note est une question posée par J. Tilouine à L. Fargues sur
l'irréductibilité des variétés d'Igusa de seconde espèce, qui lui répond que cela découle de mon
travail \cite{boy-duke}. Malheureusement l'énoncé ne se trouve pas tel quel dans loc. cit. d'où
cette note qui ne contient strictement rien de nouveau mais qui rassemble les résultats
nécessaires à la preuve.

\bigskip


\mainmatter

\renewcommand{\theequation}{\arabic{section}.\arabic{subsection}.\arabic{smfthm}}

\section{Rappels sur les variétés de Shimura associés à certains groupes unitaires et sur les variétés d'Igusa associées}

\subsection{Variétés de Shimura simples}
\label{rappels-globaux}

On considère les variétés de Shimura simples associées à des groupes unitaires telles qu'elles
sont définies dans \cite{h-t}. On reprend les notations de \cite{boy-duke}. Soit $F=F^+ E$ un corps CM, $E/\Qm$ quadratique imaginaire pure,
dont on fixe un plongement réel $\tau:F^+ \hookrightarrow \Rm$. Le groupe unitaire $G_\tau$ est alors tel que

- $G_\tau(\Rm) \simeq U(1,d-1) \times U(0,d)^{r-1}$;

- $G_\tau(\Qm_p) \simeq (\Qm_p)^\times \times \prod_{i=1}^r (B_{v_i}^{op})^\times$ où $v=v_1,v_2,\cdots,v_r$ sont les places de $F$
au dessus de la place $u$ de $E$ telle que $p=u \lexp c u$ et où $B$ est une algèbre à division centrale sur $F$ de dimension $d^2$ vérifiant certaines
propriétés, cf. \cite{h-t}, dont en particulier d'etre soit décomposée soit une algèbre à division en toute place et décomposée à la pace $v$.

\marque Pour tout sous-groupe compact $U^p$ de $G_\tau(\Am^{\oo,p})$ et $m=(m_1,\cdots,m_r) \in \Zm_{\geq
0}^r$. On pose
$$U^p(m)=U^p \times \Zm_p^\times \times \prod_{i=1}^r \ker ( \OC_{B_{v_i}}^\times \longto
(\OC_{B_{v_i}}/v_i^{m_i})^\times )$$
Pour $U^p$ assez petit on dispose d'un schéma projectif sur $\spec \OC_v$ tel que quand $U$ varie, les $X_{U^p(m)}$ forment un système projectif
dont les morphismes de transition sont finis et plats: quand $m_1=m_1'$ ils
sont en plus étales. Le système projectif $(X_{U^p(m)})_{U^p,m}$ est naturellement muni
d'une action de $G_\tau(\Am^\oo)$.

Étant donné une représentation irréductible $\xi$ de $G$ sur $\bar \Qm_l$, on notera $\LC_\xi$ le système
local sur $X_{U^p(m)}$.

\marque On note $\bar X_{U^p,m}$ la fibre spéciale de $X_{U^p(m)}$. Pour tout $0 \leq h \leq d-1$, on dispose
d'une strate fermée (resp. ouverte) notée
$$\bar X_{U^p,m}^{\geq d-h}=\bar X_{U^p,m}^{[h]} \qquad (\hbox{resp. }\bar X_{U^p,m}^{=d-h}=\bar X_{U^p,m}^{(h)})$$
de pure dimension $h$ et munie d'une action de
$G(\Am^\oo)$. Dans le cas de bonne réduction, i.e. $m_1=0$, elles sont en outre lisses. Pour tout $0<h< d$,
les strates $\bar X_{U^p,m}^{(h)}$ sont géométriquement induites sous l'action du parabolique
$P_{h,d}(\OC_v)$, au sens où il existe un sous-schéma fermé $\bar X_{U^p,m,0}^{(h)}$ tel que
$$\bar X_{U^p,m}^{(h)} \simeq \bar X_{U^p,m,0}^{(h)} \times_{P_M(\OC_v/v^{m_1})} GL_d(\OC_v/v^{m_1})$$
On notera $\bar X_{U^p,m,0}^{[h]}$ l'adhérence de $\bar X_{U^p,m,0}^{(h)}$ dans $\bar X_{U^p,m}^{[h]}$.

\subsection{Variétés d'Igusa de première espèce}

Dans \cite{h-t}, les auteurs définissent les variétés d'Igusa de première espèce
$I_{U^p,m}^{(h)} \longto \bar X_{U^p,\bar m}^{(h)}$, où $\bar m=(0,m_2,\cdots,m_r)$, qui est l'espace des
modules des isomorphismes
$$\alpha_1^{et}:(\varpi^{-m_1} \OC_v/\OC_v)^h \times \bar X_{U^p,\bar m}^{(h)} \simeq \GC^{et}[\varpi^{m_1}]$$
On obtient ainsi un revêtement galoisien de $\bar X_{U^p,\bar m}^{(h)}$ de groupe de Galois $GL_h(\OC_v/\varpi^{m_1})$. Ils définissent en outre un morphisme radiciel de $I_{U^p,m}^{(h)}$ sur
$\bar X_{U^p,m,M}^{(h)}$ pour tout $M$.

Ces variétés $I_{U^p,m}^{(h)}$ sont munies d'une action de
$$G_h^+(\Am^\oo)=G_\tau(\Am^{\oo,p} \times \Qm_p^\times \times P_{d-h,d}(F_v) \times \coprod_{i=2} (B_{v_i}^{op})^\times$$
qui agit via le quotient $G_h^+(\Am^\oo) \longto G_\tau(\Am^{\oo,p}) \times \Qm_p^\times \times
(\Zm \times GL_h(F_v))^+ \times \prod_{i=2}^r (B_{v_i}^{op})^\times$ où $(\Zm \times GL_h(F_v))^+$
désigne le sous-semi-groupe de $\Zm \times GL_h(F_v)$ formé des éléments $(c,g)$ tels que
$\varpi_o^{\lfloor - c/(n-h) \rfloor} g \in GL_h(\OC_v)$ qui est induit par la surjection
$P_{d-h,d}(F_v) \longto \Zm \times GL_h(F_v)$ qui à $g_v=g_M(g_v^c,g_v^{et})g_M^{-1}$ associe
$(v(\det g_v^c),g_v^{et})$.

\begin{theo} \label{theo1}
Pour tout $0 < h < d$, $I^{(h)}_{U^p,m}$ est géométriquement connexe.
\end{theo}

\subsection{Variétés d'Igusa de seconde espèce}

Soit $S$ un schéma et $H$ un $\OC_v$-module de Barsotti-Tate compatible de dimension $1$ et de hauteur constante $h$.
On considère comme dans \cite{h-t}, le foncteur en $S$-schémas qui à $T/S$ associe l'ensemble
des $T$-isomorphismes
$$\Sigma_{h}[\varpi^m] \times_{\spec \bar \Fm_p} T \longto H[\varpi^m] \times_S T$$
où $\Sigma_h$ est "le" $\OC_v$-module de Barsotti-Tate formel de hauteur $h$ sur $\bar \Fm_p$.
Ce foncteur est représenté par un $S$-schéma $X_m(H/S)$ de type fini sur $S$ et on note
$Y_m(H/S)$ l'intersection des images schématiques de
$$X_{m'}(H/S) \longto X_m(H/S)$$
pour $m' \geq m$ et on pose $J^{(m)}(H/S)=Y_m(H/S)^{red}$.

Soit $\GC^0/I_{U^p,m}^{(h)}$ le $\OC_v$-module de Barsotti-Tate formel universel:
il est de dimension $1$ et de hauteur constante $h$. On définit alors
$$J^{(h)}_{U^p,m,s}:=J^{(s)}(\GC^0/I_{U^p,m}^{(h)})$$
Le revêtement $J^{(h)}_{U^p,m,s} \longto I^{(h)}_{U^p,m} \times \spec \bar \Fm_p$ est galoisien de groupe de
Galois $(\DC_{v,d-h}/ \varpi^s)^\times$ où $\DC_{v,d-h}$ est l'ordre maximal de
l'algèbre à division $D_{v,d-h}$ de centre $F_v$ et d'invariant $1/(d-h)$.

Par ailleurs, on a, sur la tour des $J^{(h)}_{U^p,m,s}$ une action par correspondances du groupe
$G_\tau(\Am^{\oo,p}) \times \Qm_p^\times \times (\Zm \times GL_h(F_v))^+ \times \prod_{i=2}^r
(B_{v_i}^{op})^\times \times (D_{v,d-h}^\times/\DC_{v,d-h}^\times)$.

\begin{defi} Soit $\DC_{v,d-h}^1$ le noyau de la norme réduite de $D_{v,d-h}^\times$. Pour tout $m \geq 0$,
on note alors $\DC_{v,d-h,m}^1:=\DC_{v,d-h}^1/(1+\varpi^s)$ et soit $J^{(h),1}_{U^p,m,s}:=
(J^{(h)}_{U^p,m})^{\DC_{v,d-h,m}^1}$.
\end{defi}

\begin{theo} \label{theo2}
Pour tout $0 < h < d$, $J^{(h)}_{U^p,m} \longrightarrow J^{(h),1}_{U^p,m}$ est un
$\DC_{v,d-h,m}^1$-torseur géométriquement irréductible au dessus de toute composante
géométriquement connexe et
$$J^{(h),1}_{U^p,m}(\bar \Fm_p) \simeq I^{(h)}_{U^p,m}(\bar \Fm_p) \times \DC_{v,d-h,m}/ \DC_{v,d-h,m}^1$$
où l'action est induite par l'action du groupe de Galois $\DC_{v,d-h,m}$ du revêtement
$J^{(h)}_{U^p,m} \longto I^{(h)}_{U^p,\bar m}$.
\end{theo}

\section{Systèmes locaux d'Harris-Taylor: groupes de cohomologie}

\subsection{Définition}

Les variété d'Igusa de seconde espèce définissent pour toute représentation irréductible
admissible $\rho_v$ des inversibles $D_{v,d-h}^\times$, un système local $\FC_{\rho_v,U^p(m),0}$
sur $I_{U^p,\bar m}^{(h)}$ et donc sur $\bar X_{U^p,m,0}^{(d-h)}$. Tout élément
$(g^p,g_{p,0},c,g_o^{et},g_{o_i},\d)$ de $G_\tau(\Am^{\oo,p}) \times \Qm_p^\times \times (\Zm
\times GL_{d-h}(F_v))^+ \times \prod_{i=2}^r (B_{v_i}^{op})^\times \times
(D_{v,d-h}^\times/\DC_{v,d-h}^\times)$ définit naturellement un morphisme
$$(g^p,g_{p,0},c,g_v^{et},g_{v_i},\d):(g^p,g_{p,0},c,g_v^{et},g_{v_i},\d)^* (\FC_{\rho,M} \otimes \LC_\xi)
\longto \FC_{\rho,M} \otimes \LC_\xi$$

\begin{defi}
Pour $\tau_v$ une $\bar \Qm_l$-représentation irréductible admissible de $D_{v,h}^\times$, on note
$$\FC_{\tau_v}:= \FC_{\tau_v,0} \times_{P_{h,d}(\OC_v/v^{m_1})} GL_d(\OC_v/v^{m_1})$$
le faisceau induit sur toute la strate $\bar X_{U^p,m}^{(d-h)}$, où l'action de $\left (
\begin{array}{cc} g_v^c & * \\ 0 & g_v^{et} \end{array} \right) \in P_{h,d}(F_v)$
est donnée par celle de $(\val (\det g_v^c),g_v^{et}) \in \Zm \times GL_{d-h}(F_v)$.
\end{defi}

\subsection{Faisceaux pervers d'Harris-Taylor}

On introduit suivant \cite{boy-duke}, les faisceaux pervers dits de Harris-Taylor
$$\PC(g,t,\pi_v)=j^{\geq tg}_{!*} HT(g,t,\pi_v,[\overleftarrow{t-1}]_{\pi_v})[d-tg]
\otimes \rec_{F_v}^\vee(\pi_v)$$ où $HT(g,t,\pi_v,\Pi)=
\FC_{\JL^{-1}([\overleftarrow{t-1}]_{\pi_v})} \otimes \Pi$. On rappelle le résultat suivant de
\cite{boy-duke}.

\begin{prop} Pour tout $1 \leq t \leq s$, la filtration par les poids de $j^{\geq tg}_! HT(g,t,\pi_v,
\Pi_t)[d-tg]$ est telle que ses gradués $gr_{!_t,k}$ sont nuls pour $k < 0$ et pour $0 \leq k \leq s-t$, égaux à
$$j^{\geq (t+k)g}_{!*} HT(g,t+k,\pi_v,\Pi_t \overrightarrow{\times} [\overleftarrow{k-1}]_{\pi_v})[d-(t+k)g]
\otimes \Xi^{k(g-1)/2}$$
où $\Xi:\Zm \longto \bar \Qm_l^\times$ désigne le caractère défini par $\Xi(1)=\frac{1}{l}$

\end{prop}

\subsection{Preuve des théorèmes (\ref{theo1}) et (\ref{theo2})}

Pour tout $\tau_v$ et pour tout $i$,
$$H^i_{\tau_v,0}:=H^i(\lim_{\atop{\longto}{U^p,m}} \bar X_{U^p,m,0}^{(d-h)},
\FC_{\tau_v,0})$$ est une représentation admissible de $G(\Am^{\oo,o}) \times GL_{h}(F_v) \times
\Zm$ telle que pour toute représentation $\Pi_h$ de $GL_{d-h}(F_v)$, on a
$$H^i_{\tau_v}(\Pi):=H^i(\lim_{\atop{\longto}{U^p,m}} \bar X_{U^p,m}^{(d-h)}, \FC_{\tau_v}\otimes \Pi_h) \simeq
\Ind_{P_{h,d}(F_v)}^{GL_d(F_v)} (H^i_{\tau_v,0} \circledast \Pi_h)$$ où l'action de $\left (
\begin{array}{cc} g_v^c & * \\ 0 & g_v^{et} \end{array} \right) \in P_{h,d}(F_v)$ sur
$H^i_{\tau_v,0} \circledast \Pi$ est donnée par l'action de $(\val (\det g_v^c),g_v^{et}) \in \Zm
\times GL_{d-h}(F_v)$ sur $H^i_{\tau_v,0}$ et par celle de $g_v^c$ sur $\Pi_h$. Les théorèmes
(\ref{theo1}) et (\ref{theo2}) découlent classiquement du résultat suivant.

\begin{prop} Pour tout $\tau_v$, $H^0_{\tau_v,0}$
est nul sauf si $\tau_v$ est un caractère $\chi_v$ auquel cas on obtient
$$\chi_v \otimes \Xi^0$$
en tant que représentation de $GL_{d-h}(F_v) \times \Zm$.
\end{prop}

\begin{proof}
On rappelle qu'une représentation automorphe $\Pi$ de $G_\tau(\Am)$ est dite cohomologique, s'il existe une certaine représentation algébrique $\xi$ sur $\Cm$ de la restriction des scalaires de $F$ à $\Qm$ de $GL_g$, et un entier $i$ tels que
$$H^i((\lie G_\tau(\Rm)) \otimes_\Rm \Cm,U_\tau,\Pi_\oo \otimes (\xi')^\vee) \neq (0)$$
où $U_\tau$ est un sous-groupe compact modulo le centre de $G_\tau(\Rm)$, maximal, cf. \cite{h-t} p.92, et où $\xi'$
est le caractère sur $\Cm$ associé à $\xi$ via un isomorphisme $\bar \Qm_l \simeq \Cm$ fixé.
Considérons $\Pi$ irréductible et cohomologique, les faits suivants se trouvent dans \cite{boy-duke}:

\begin{itemize}
\item si $\Pi$ n'est pas cohomologique alors $H^i_{\tau_v,0}[\Pi^{\oo,v}]$ est nul pour tout $i$ et tout $\tau_v$;

\item pour $\Pi$ cohomologique tel que $\Pi_v \simeq \st_s(\pi_v)$, pour $d=sg$ avec $\pi_v$ une représentation
irréductible cuspidale de $GL_g(F_v)$, $H^i_{\JL^{-1}(\st_t(\pi_v)),0}$ est donné par
$$\left \{ \begin{array}{cl} 0 & \hbox{si } i \neq 0 \\
\sharp \ker^1(\Qm,G_\tau) m(\Pi) [\overleftarrow{s-t-1}]_{\pi_v(-t(g-1)/2)} \otimes ( \Xi^{\frac{(s-t)(g-1)}{2}}
\otimes \bigoplus_{\chi \in \AF(\pi_v)} \xi^{-1}) & i=0 \end{array} \right.$$
où $\AF(\pi_v)$ est l'ensemble des caractères
$\chi:\Zm \longto \Qm_l^\times$, tels que $\pi_v \otimes \chi^{-1} \circ \val(\det) \simeq \pi_v$.

\item pour $\Pi$ cohomologique tel que $\Pi_v \simeq \st_s(\pi_v)$, pour $d=sg$ avec $\pi_v$ une représentation
irréductible cuspidale de $GL_g(F_v)$, $H^i_{\JL^{-1}(\st_t(\pi_v)),0}$ est donné par
$$\left \{ \begin{array}{cl} 0 & \hbox{si } i \neq s-t \\
\sharp \ker^1(\Qm,G_\tau) m(\Pi) [\overrightarrow{s-t-1}]_{\pi_v(-t(g+1)/2} \otimes (\Xi^{\frac{(s-t)(g+1)}{2}}
\bigoplus_{\chi \in \AF(\pi_v)} \chi^{-1}) & i=s-t \end{array} \right.$$

\item pour $\Pi$ cohomologique telle que $\Pi_v$ n'est pas d'une des deux formes précédentes,
d'après loc. cit., $\Pi_v$ est alors de la forme
$$[\overrightarrow{s_1-1}]_{\pi_{v,1}} \times \cdots \times [\overrightarrow{s_u-1}]_{\pi_{v,u}} \times
[\overleftarrow{t_1-1}]_{\pi_{v,1}'} \times \cdots \times
[\overleftarrow{t_{u'}-1}]_{\pi_{v,u'}'}$$ avec $\pi_{v,i}$ et $\pi_{v,j}'$ irréductibles
cuspidales unitaires de respectivement $GL_{g_i}(F_v)$ et $GL{g_j'}(F_v)$. \footnote{En outre tous
les $s_i$ ont la même congruence modulo $2$.} On a alors prouvé dans loc. cit. que
$$H^0(\lim_{U^p,m} \bar X_{U^p,m}^{[d-h]}, \PC(g,t,\pi_v))[\Pi^{\oo,v}]=0.$$
En utilisant la suite spectrale associée à la filtration par les poids de $j_!^{\geq
tg}HT(g,t,\pi_v,[\overleftarrow{t-1}]_{\pi_v})[d-tg]$, on en déduit que
$$H^0_{\tau_v}([\overleftarrow{t-1}]_{\pi_v})[\Pi^{\oo,v}]=0$$
et donc que $H^0_{\tau_v,0}[\Pi^{\oo,v}]=0$, d'où le résultat.

\end{itemize}

\end{proof}

\bibliographystyle{plain}
\bibliography{bib-ok}

\end{document}